\documentclass[11pt]{amsart}
\usepackage{amsmath,amsfonts,amssymb,amsthm,epsfig,color}
\newcommand{\bb}{\mathbb}

\newcommand{\h}{\bb H}
\newcommand{\Z}{\bb Z}
\newcommand{\R}{\bb R}
\newcommand{\N}{\bb N}

\newcommand{\Q}{\bb Q}

\newcommand{\f}{\mathfrak}
\newcommand{\M}{\mathcal M}

\newcommand{\T}{\mathcal T}
\newcommand{\hh}{\mathcal H}

\newtheorem{Theorem}[equation]{Theorem}
\newtheorem{Cor}[equation]{Corollary}
\newtheorem{Prop}[equation]{Proposition}
\newtheorem{lemma}[equation]{Lemma}
\newtheorem*{lemma*}{Lemma}

\newtheorem*{theorem*}{Theorem}
\newtheorem{Def}[equation]{Definition}
\numberwithin{equation}{section}
\begin{document}
\title[Logarithm laws]{Logarithm laws and shrinking target properties}
\author{J.~S.~Athreya}
\email{jathreya@math.princeton.edu}
\address{Dept. of Mathematics, Princeton University, Princeton, NJ}
\thanks{J.S.A. supported by NSF grant DMS
    0603636.}

\begin{abstract}
We survey some of the recent developments in the study of logarithm laws and shrinking target properties for various families of dynamical systems. We discuss connections to geometry, diophantine approximation, and probability theory.
\end{abstract}
\maketitle

\begin{center}\textit{Dedicated to Professor S.~G.~Dani on his 60th birthday}\end{center}

\section{Introduction}

A common theme in the study of (deterministic) dynamical systems is to ask when their long-term behavior mimic that of random processes. For example, the Birkhoff ergodic theorem and its generalizations are `laws of large numbers', and mixing properties can be viewed as questions about `long-term independence'. 

A natural question in this spirit is to ask when a system satisfies a `Borel-Cantelli'-type property. Recall that the classical Borel-Cantelli lemma (see section~\ref{shift} for a precise statement) states that a sequence of independent events occurs infinitely often with probability 1 if the sum of their probabilities diverges, and occurs finitely often with probability one if the sum of their probabilities converge. The second statement does not require independence, and thus can be immediately extended to deterministic systems. 

We will consider the following general setup: let $(X,  \mathcal{F}, \mu)$ be a probability space, and $G$ be a group acting on $X$ by measure-preserving transformations. Let $\mathcal{A} \subset \mathcal{F}$ be a family of measurable subsets of $X$.  

\begin{Def}\label{STP} Let $\{g_n\}_{n \in \N} \subset G$ be a sequence of group elements. We say that a sequence of measurable subsets $\{A_n\}_{n \in \N}$ is a \emph{Borel-Cantelli (BC) sequence} for $\{g_n\}$ if $\sum_{n=1}^{\infty} \mu(A_n) = \infty$ and $$\mu(\{x \in X : g_n x \in A_n \mbox{ infinitely often} \}) =1.$$We say it is a \emph{strongly Borel-Cantelli (sBC) sequence} if $$\lim_{n \rightarrow \infty}\frac{ |\{1 \le i \le n: g_i x \in A_i\}|}{\sum_{i=1}^n \mu(A_i)} = 1$$ for almost every $x \in X$. We say a family of measurable subsets $\mathcal{A}$ is \emph{(strongly) Borel-Cantelli ([s]BC)} for $\{g_n\}$ if for every sequence of subsets $\{A_n\}_{n \in \N} \subset \mathcal{A}$ with divergent measure is a (strongly) Borel-Cantelli sequence. \end{Def}

Note that if $\mathcal{A}$ is Borel-Cantelli for $\{g_n\}$, we have, for any sequence $\{A_n\}_{n \in \N} \subset \mathcal{A}$ \begin{equation}
\mu(\{x \in X : g_n x \in A_n \mbox{ infinitely often} \}) =
\left\{ \begin{array}{ll}  1 & \sum_{n=0}^{\infty} \mu(A_n) = \infty \\ 0 & \mbox{otherwise}\end{array}\right.\nonumber\end{equation}

The convergence half of this equation follows from the (easy half) of the classical Borel-Cantelli Lemma~\ref{bc}. Often, we are interested in the situation when our target sets $\{A_n\}_{n \in \N}$ are decreasing, or ``shrinking'' in an appropriate manner: hence the phrase ``shrinking target property''. 

This terminology was first coined by Hill-Velani~\cite{HilVel}, who studied the Hausdorff dimensions of the $\limsup$ set $\{x \in X : g_n x \in A_n \mbox{ infinitely often} \}$ in the \emph{convergence} cases, when $\sum_{n=0}^{\infty} \mu(A_n) < \infty$ (and thus this set has measure $0$).

We will be primarily interested in the \emph{divergence} case in two main contexts:

\subsection{$\Z$-actions on metric spaces} Here, we are interested in the action of a transformation $T$ preserving a Borel probability measure $\mu$ on a metric space $X$, and we will primarily focus on metric balls as our target sets. A classical example of a theorem in this context, due to Phillip~\cite{Phillip}, is one of the first Borel-Cantelli laws in a deterministic context:

\begin{Theorem}\label{Phillips} If $T: [0, 1) \rightarrow [0, 1)$ is either the expanding map $T(x) = \beta x \mbox{ (mod } 1)$ ($\beta >1$) or the Gauss map $T(x) = \{1/x\}$. Then the family of intervals is sBC.
\end{Theorem}

In general, following~\cite{Fayad}, we say a transformation $T$ of a metric space $X$ has the \emph{shrinking target property} (STP) if the family of metric balls is BC, and that it has the \emph{monotone shrinking target property} (MSTP) if every sequence of balls with decreasing radii and divergent sum of measure is a BC sequence. In this language, the above result can be rephrased as saying that the maps discussed satisfy STP.

\subsection{Homogeneous actions} Here, we consider actions of a Lie group (and its various subgroups) $G$ on a homogeneous space $X = G/\Gamma$, where $\Gamma \subset G$ is a lattice. 
$G$ acts by translations on $X$, preserving the (normalized) Haar measure $\mu$. 

In the setting where $\Gamma$ is non-uniform, it is natural to study the family of ``cusp neighborhoods'' given by \emph{complements} of balls. For example, consider the case where $G = SL(2, \R)$, and $(X, \mu)$ is the unit-tangent bundle to the non-compact, finite volume hyperbolic surface $\h^2/\Gamma$ with Lebesgue measure $\mu$. Here, the actions of the subgroups

 $$A = \left\{g_t = \left(\begin{array}{cc} e^{t/2} & 0 \\ 0 & e^{-t/2}
\end{array}\right)\right\}_{t \in \R};$$\noindent and $$U = \left\{h_t =\left(\begin{array}{cc} 1 & t \\ 0 & 1
\end{array}\right)\right\}_{t \in \R},$$

\noindent correspond to the geodesic and horocycle flow respectively. Each flow is ergodic, and thus a generic orbit is dense (in fact, for the horocycle flow the results of Dani-Smillie~\cite{DS} tell us that all non-periodic orbits are dense).

The statistical behavior of their excursions into the cusp can be described by the following theorem, known as a \emph{logarithm law}, due in the case of geodesic flow to Sullivan~\cite{Sullivan} and for horocycle flow to Athreya-Margulis~\cite{AM}. 

\begin{Theorem}\label{loglawsl2} Let $G = SL(2, \R)$, $\Gamma \subset G$ a non-uniform lattice. Let $\{\varphi_t\}_{t \in \R}$ denote either horocycle or geodesic flow. Then for all $y \in X$, and $\mu$-almost every $x \in X$, \begin{equation}\label{loglaw}\limsup_{t \rightarrow \infty} \frac{ d(\varphi_t(x),y)}{\log t} = 1,\end{equation} where distance is measured in the projection to the hyperbolic surface $\h^2/\Gamma$.    \end{Theorem}

\subsection{Organization}In the rest of this survey, we try and give our perspective the recent flowering of results in the area of logarithm laws and shrinking target properties, much of which has been inspired by the papers~\cite{KM, Sullivan}. In section~\ref{hyp}, we describe these results for geodesic flows on homogeneous spaces, spaces of negative curvature, and Teichmuller space. In section~\ref{shift}, we study such theorems for (partially) hyperbolic diffeomorphisms and shift maps, as well as recaling the classical Borel-Cantelli lemma and some of its generalizations. In sections~\ref{ellip} and~\ref{para}, we describe some results in the elliptic and parabolic contexts, respectively. Here, in lieu of rapid mixing, diophantine properties and the geometry of numbers come into play. Finally, in section~\ref{conj}, we discuss some of the open problems and conjectures in the area.

\subsection{Acknowledgments} I thank Krishna Athreya, Steven Lalley, Gregory Margulis, Yair Minsky, and Nimish Shah for useful discussions; and Anish Ghosh and Dmitry Kleinbock for their detailed comments. David Ellwood and the Clay Mathematics Institute have provided both financial and organizational support. I also thank the Tata Institute of Fundamental Research, and in particular the organizers of the International Conference on Dynamics and Measures on Groups and Homogeneous Spaces for their hospitality.

Finally, I would like to take this opportunity to thank and personally acknowledge Professor S.~G.~Dani. The results in this paper are at the nexus of geometry, dynamics, probability theory, and diophantine approximation; and Professor Dani has been at the forefront of exploring these connections.  On a more personal note, his lectures in the Fall of 2000 at the University of Chicago, when I was in my first year of graduate school, were influential in my choice of area of study.

\section{Geodesic flows}\label{hyp}

\subsection{Homogeneous spaces}

Theorem~\ref{loglawsl2} can be generalized in many directions. In fact, Sullivan's original result in~\cite{Sullivan} is for the geodesic flow for a general non-compact finite volume hyperbolic manifold.

\begin{Theorem}\label{sullivan} Let $\Gamma \subset SO(k+1,1)$ be a non-uniform lattice, $x, y \in X = \h^{k+1}/\Gamma$.  Let $\{g_t\}$ denote the one parameter subgroup corresponding to geodesic flow on the unit tangent bundle $SX$, and $p: SX \rightarrow X$ the projection map. Then for almost all (with respect to Lebesgue measure on the sphere) $\theta \in S_x X$, $$\limsup_{t \rightarrow \infty} \frac{d(p(g_t(x, \theta)), y)}{\log t} = \frac{1}{k}.$$
\end{Theorem}

This in fact a corollary of a Borel-Cantelli result: the family of sets given by complements of balls is strongly Borel-Cantelli for the sequence $\{g_n\}$.

This paper was the first (to this author's knowledge) to use the terminology ``logarithm law." Sullivan's proof is quite subtle, and relies on a Borel-Cantelli type theorem for spheres in Euclidean space- actually horoballs corresponding to cusp neighborhoods in $\h^{k+1}/\Gamma$. 

For $k=1$, and $\Gamma = SL(2, \Z)$, these horoballs are simply the interior of the \emph{Ford circles}. Recall that the Ford circle based at a rational number $p/q \in \Q$ is the circle tangent to the real axis at $p/q$ with radius $1/2q^2$. In fact, in this case, the Sullivan result can be derived as a corollary to the classical Khinchin theorem on metric diophantine approximation. 

\vspace{.1in}

\noindent \textbf{Remark:} There is in fact a natural correspondence between logarithm laws, shrinking target properties, and diophantine approximation. There is a rich literature and much active research exploring these connections which is not in the scope of this paper. For a very interesting recent example of how questions on diophantine approximation can be reduced to studying shrinking target properties for flows on homogeneous spaces, see the paper~\cite{GS} of Gorodnik-Shah.

\vspace{.1in}

Motivated by a diophantine connection, Kleinbock-Margulis~\cite{KM} extended the result of Sullivan to non-compact Riemannian symmetric spaces: 

\begin{Theorem}\label{loglawzero}(\cite{KM}, Theorem 1.7 and Prop 5.1)  Let $G$ be a connected, center-free, semisimple Lie group without compact factors, $\f{g}$ its Lie algebra, $\Gamma \subset G$ an irreducible non-uniform lattice, $K$ a maximal compact subgroup, and $d(.,.)$ a distance function on $G/\Gamma$ determined by a right-invariant Riemannian metric on $G$ bi-invariant under $K$. Let $\mu$ denote the measure on $G/\Gamma$ determined by Haar measure on $G$. Let $\f{a} \subset \f{g}$ be  a Cartan subalgebra, $z \in \f{a}$, and $g_t = \exp(tz)$. Then there exists a $k = k(G/\Gamma) \in (0, \infty)$ such that $\forall y$, for $\mu$-a.e. $x$, $$\limsup_{t \rightarrow \infty} \frac{d(g_t x, y)}{\log t} = 1/k$$\end{Theorem}

\noindent  Again, this is a corollary of a shrinking target result: the family of sets given by complements of balls is strongly Borel-Cantelli for the sequence $\{g^n\}$.

\noindent Here, $g = g_1 = \exp(z)$. The proof of this result has two main ingredients: 

\begin{itemize} 

\item Mixing properties of the $G$-action on $G/\Gamma$, in particular, quantitative estimates on the decay of matrix coefficients following the Howe-Moore theorem. 

\item An estimate on the (exponential) rate of decay of the volumes of cusp neighborhoods $$\mu(x: d(x, y) > t).$$

\end{itemize}

These are provided by:

\begin{Theorem}\label{matrix}(\cite{KM}, Theorem 1.12) Let $\rho_0$ denote the regular representation of $G$ on the subspace of $L^2(G/\Gamma)$ orthogonal to the constant functions. Then the restriction of $\rho_0$ to any simple factor of $G$ is isolated (in the Fell topology) from the trivial representation.\end{Theorem}

\begin{Prop}\label{cusp}(\cite{KM}, Proposition 5.1) For all $y \in G/\Gamma$, there are constants $C_1, C_2 >0$ such that $$C_1 e^{-kt} \le \mu(x \in G/\Gamma: d(x, y) > t) \le C_2 e^{-kt}.$$ \end{Prop}

The idea is to use Theorem~\ref{matrix} to obtain mixing estimates to say that the pre-images of the sets $$A_n : = \{x: d(g_n x, y) > r_n\}$$\noindent become almost independent very rapidly. Then, one can apply the Borel-Cantelli lemma (Lemma~\ref{bc}), or more precisely an appropriate generalization (Lemma~\ref{sprindzuk}). Then the problem reduces to estimating the measures of the sets $A_n$, which we can do by Proposition~\ref{cusp}.

Proposition~\ref{cusp} is proved by studying the reduction theory of algebraic groups, and in many cases the constant $k$ can be explicitly calculated using systems of roots.

\subsection{Trees, buildings, and positive characteristic}\label{trees}

Let $k$ be a local field of characteristic $p>2$, and $G$ be the $k$-points of a semisimple linear algebraic group defined over $k$, and $\Gamma \subset G$ a non-uniform lattice. Let $X = X_G$ denote the Bruhat-Tits building of $G$, $\partial X$ the geodesic boundary of $X$,  $Y = X/\Gamma$, and $\pi: X \rightarrow Y$ be the natural projection. Fix a $G$-invariant metric $d$ on $X$ (we will also use $d$ to denote the metric on the quotient $Y$). Let $\{ g_t (x, \theta)\}_{t \geq 0}$ denote the geodesic starting at  the vertex $x \in X$ in direction $\theta \in \partial X$. There is a natural measure class on $\partial X$. 

\begin{Theorem}\label{loglawp}~\cite{AGP} There is a $k = k(Y)>0$ such that  for any $x \in X$, $y \in Y$ and almost all $\theta \in \partial X$, $$\limsup_{t \rightarrow \infty} \frac{ \log d_{Y} (\pi(g_t(x, \theta)), y)}{ \log n} = 1/k.$$
\end{Theorem}

Once again, the idea is to combine volume and mixing estimates: reduction theory of algebraic groups to obtain an analogue of Proposition~\ref{cusp}, and estimates on decay of matrix coefficients.

If our ambient group has $k$-rank $1$, the building is a tree. However, the algebraic group may not be the entire automorphism group of the tree, and there are interesting non-algebraic quotients. These were studied by Hersonsky-Paulin~\cite{HP}. Let $\T$ be a locally finite tree, and $G = \mbox{Aut}(T)$. Let $\Gamma \subset G$ be a non-uniform lattice. There is a natural height function $h = h_{\Gamma}: \T \rightarrow (0, \infty)$ that is essentially equivalent to the distance from a fixed point in the quotient $\T/\Gamma$. Let $\delta>0$ be the Hausdorff dimension of the boundary $\partial\T$, and $\mu$ the associated Hausdorff measure. Given a vertex $x \in \T$, and a point $\eta \in \partial\T$, let $c_{\eta}$ be the geodesic ray starting at $x$ in direction $\eta$. 

\begin{Theorem}\label{loglawtrees}~\cite{HP} For all $x \in \T$, $\mu$-almost every $\eta \in \partial \T$, $$\limsup_{t \rightarrow \infty} \frac{h(c_{\eta}(t))}{\log t} = 1/\delta.$$\end{Theorem}

\subsection{Negative curvature}\label{negcurv}


Theorem~\ref{loglawtrees} is one part of a larger program of Hersonsky-Paulin. In a remarkable series of papers~\cite{HP0, HP, HP3}, they have studied the statistical properties of geodesics in spaces of negative curvature, including, but not limited to, excursions into cusps.  

Their methods are also very different from those in the homogeneous setting. They follow more closely the original arguments of Sullivan, relying on subtle analysis of the geometry of horoballs (when studying cusps).  One of their main results is a generalization of Theorem~\ref{sullivan} to the situation of  variable negative curvature Riemannian manifolds (an earlier generalization, for $M$ geometrically finite with constant curvature, is due to Stratmann-Velani~\cite{StratVel}).

Let $M$ be a complete pinched negatively curved geometrically finite Riemannian manifold with one cusp, $e$. Given a Riemannian covering $\pi: N \rightarrow M$, and $x_0 \in N$, let $$f_{\pi}(t) = |\{x \in \pi^{-1}(\pi(x_0)): d_N(x, x_0) < t\}|$$\noindent be the associated counting function.

Let $\tilde{\pi}: \tilde{M} \rightarrow M$ be a universal covering, and $\pi_0: M_0 \rightarrow M$ be the covering of $M$ associated to a parabolic subgroup $H_e \subset \pi_1(M)$ corresponding to the cusp $e$. Suppose there are $\delta > \delta_0 >0$ such that $$f_{\tilde{\pi}}(t) \sim e^{\delta t},  f_{\pi_0}(t) \sim e^{\delta_0 t}.$$\noindent Here $f \sim g$ indicates there is a constant $C \geq 1$ such that $1/C \le f/g \le C$. $\delta$ is called the critical exponent of the manifold $M$. Given $x \in M$, let $\mu$ be the Patterson-Sullivan measure on the unit tangent sphere $S_x(M)$. 

\begin{Theorem}\label{loglawhyp}~\cite{HP0} For all $x, y \in M$, for $\mu$-almost every $\theta \in S_x M$, $$\limsup_{t \rightarrow \infty} \frac{d(g_t(x, \theta), y)}{\log t} = \frac{1}{2(\delta - \delta_0)}.$$\noindent Here $g_t(x, \theta)$ is the geodesic determined by $\theta \in S_x M$\end{Theorem}

\noindent In~\cite{HP3}, they investigate spiraling of geodesic trajectories around closed geodesics and other totally geodesic submanifolds, proving shrinking target properties and logarithm laws for neighborhoods of those objects. These are different from many of the other results described in this paper in a qualitative way, in that they study visits to a compact part of the state space, rather than a cusp or neighborhood of infinity. 

Perhaps the most natural family of compact sets are metric balls. For geodesic flow on a hyperbolic manifold (i.e., the setting of Sullivan's result), Maucourant~\cite{Mau} studied their shrinking target properties, deriving a general Borel-Cantelli lemma, which can be viewed as an MSTP (see the introduction) for geodesic flow.

\begin{Theorem}\label{mac}~\cite{Mau} Fix notation as in Theorem~\ref{sullivan}. Let $\{B_t\}_{t \geq 0}$ be a decreasing family of closed metric balls of radius $\{r_t\}_{t \geq 0}$. Then for all $x \in X$, 
\begin{equation}\label{hypballstp}\lambda\{\theta \in S_x X  :  \{ t \geq 0: \pi(g_t(x, \theta)) \in B_t\} \mbox{ is unbounded } \} =
\left\{ \begin{array}{ll}  1 & \int_{t =0}^{\infty} r_t^{k} = \infty \\ 0 & \mbox{otherwise}\end{array}\right.\end{equation}\end{Theorem}

As a corollary, one obtains a logarithm law for visits to balls: for all $x, y \in X$, almost every $\theta \in S_x X$, $$\limsup_{t \rightarrow \infty} \frac{-\log d(\pi(g_t(x, \theta)), y)}{\log t} = 1/k.$$

We will study more result of this type (i.e., on visits to compact parts of the space) in section~\ref{shift}. In contrast to the results of Chernov-Kleinbock~\cite{CK} and Dolgopyat~\cite{Dolgopyat} that we will discuss there, Maucourant does not explicitly use rapid mixing of the geodesic flow, instead proving a generalization of the Borel-Cantelli lemma based on controlling $L^p$-norms of certain sums, and using properties of hyperbolic geometry to ensure this generalization can be applied.

\subsection{Teichm\"uller geodesic flow}

While it is not symmetric, not negatively curved, and in fact not even Riemannian, the moduli space of Riemann surfaces endowed with the Teichm\"uller metric does share many important geometric and dynamical properties with non-compact, finite volume hyperbolic manifolds and symmetric spaces. One such property is the logarithm law for geodesic flow, due in this setting to Masur~\cite{Masur}. 

Let $\M_g$ denote the moduli space of compact genus $g$ Riemann surfaces, and $Q_g$ denote the moduli space of holomorphic quadratic differentials. This is an orbifold vector bundle over $\M_g$, in fact, it is the cotangent bundle. $Q_g^{(1)} \subset Q_g$, the sub-bundle of unit-area differentials, corresponds to the unit cotangent bundle. Let $\pi: Q_g^{(1)} \rightarrow \M_g$ denote the natural projection. There is a natural $SL(2,\R)$-action on $Q_g^{(1)}$. Given $q \in \pi^{-1}(X)$, the orbit of the compact subgroup $$K=\left\{r_{\theta} = \left( \begin{array}{cc} \cos \theta & \sin \theta \\
-\sin \theta & \cos \theta \end{array}\right): 0 \le \theta <
2\pi\right\}$$ \noindent yields a family $\{q_{\theta}\} \subset \pi^{-1}(X)$. Let $\{X_{\theta, t}\}_{t \geq 0}$ be the Teichm\"uller geodesic determined by $q_{\theta}$. 

\begin{Theorem}\label{teichloglaw}~\cite{Masur} Let $q \in Q_g^{(1)}$, $X = \pi(q)$. For Lebesgue almost every $\theta \in [0, 2\pi)$, $$\limsup_{t \rightarrow \infty} \frac{d(X_{\theta, t}, X)}{\log t} = 1,$$\noindent where $d$ denotes the Teichm\"uller metric on $\M_g$.\end{Theorem}

Masur's techniques are also geometric in nature, based on an analysis of the length of curves on the family of surfaces $X_t$ and their relationship to the Teichm\"uller distance.

\section{Abstract BC properties}\label{shift}

\noindent Let $(X, \mathcal{F}, \mu)$ be a probability space. Let $X_n: X \rightarrow \{0,1\}$ be a sequence of $0-1$ random variables on $S$, with $\mu(x: X_n(x) = 1) =: p_n$. Also define $p_{ij} := \mu(x: X_i(x)X_j(x) =1)$. Given measurable $f: X \rightarrow \R$, we write $E(f) := \int_{X} f d\mu$ for the expectation. Recall the classical Borel-Cantelli lemma:

\begin{lemma}\label{bc}(Borel-Cantelli) 
\begin{enumerate}
\item If $\sum_{n=0}^{\infty} p_n < \infty$, $\mu(\sum_{n=0}^{\infty} X_n = \infty) =0$
\item If the $X_n$'s are pairwise independent (i.e. $p_{nm}  = p_n p_m \forall m, n$), and $\sum_{n=0}^{\infty} p_n = \infty$,  then $\mu(\sum_{n=0}^{\infty} X_n = \infty) =1$
\end{enumerate}
\end{lemma}

The first example of a logarithm law can be derived from this as follows: Let $\{Y_n\}_{n=0}^{\infty}$'s be i.i.d. exponential $\lambda$ random variables ($\lambda \in \R^{+}$ a fixed parameter), i.e., $\mu(Y_n > t) = e^{-\lambda t}$ for all $t \geq 0$. Let $\{r_n\}_{n=0}^{\infty}$ be a sequence of positive real numbers. Applying the above results to the sequence of random variables  \begin{equation}X_n : = \left\{
\begin{array}{ll}  1 & Y_n > r_n\\ 0 & \mbox{otherwise}\end{array}\right.\nonumber\end{equation} Thus, $Y_n > r_n$ infinitely often with probability $1$ if and only if $\sum_{n=0}^{\infty} e^{-\lambda r_n} = \infty$. As a corollary, one obtains that, with probability $1$,  $$\limsup_{n \rightarrow \infty} \frac{Y_n}{\log n} = 1/\lambda.$$

More generally, the convergence result yields the upper bounds in the logarithm laws. To obtain the lower bounds/infinite cases, we have to work to relax the independence assumption in the divergence result. Let $$S_n : = \sum_{i=1}^{n} X_i \mbox{ ; }E_n : =  E(S_n) = \sum_{i=1}^n p_i.$$ We will use the following generalization of Lemma~\ref{bc}, which can be found in ~\cite{KM}:

\begin{lemma}\label{sprindzuk}(Sprindzuk) Suppose there is a constant $C$ such that for all $M > N \geq 1$, we have $$\sum_{i, j=N}^{M} |p_{ij} - p_i p_j| \le C \sum_{i=N}^{M} p_i.$$ Then for any $\epsilon >0$ we have, as $N \rightarrow \infty$, $$S_N = E_N + O(E_N^{1/2} \log^{3/2 + \epsilon}E_N),$$ with probability $1$. In particular, this yields that $S_N/E_N \rightarrow 1$ with probability $1$ if $E_N \rightarrow \infty$.\end{lemma}

Specializing definition~\ref{STP} to $\Z$-actions, we consider invertible measure-preserving maps $T: (X, \mu) \rightarrow (X, \mu)$, a family of measurable subsets $\mathcal{A}$, and given a sequence of measurable subsets $\{A_n\}_{n \in \N} \subset \mathcal{A}$, we define random variables $X_n (x) = \chi_{A_n}(T^n x)$. 

The study of the Borel-Cantelli property reduces to studying the behavior of the associated random variables $S_n$ and expectations $E_n$. Notice that the Sprindzuk lemma above yields stronger information than the fact that $S_n \rightarrow \infty$ with probability one if $E_n \rightarrow \infty$: it yields, for example, that $S_n/E_n \rightarrow 1$. We say that a sequence of sets $A_n$ is \emph{strongly Borel-Cantelli} (abbreviated sBC) if the ratio $S_n/E_n \rightarrow 1$ with probability $1$.

One can characterize many different ergodic properties of the system $T: X \rightarrow X$ in terms of the BC and sBC properties. We recall Proposition 1.5 of~\cite{CK}:

\begin{Theorem}\label{CKprop}~\cite{CK} Let $T: (X, \mu) \rightarrow (X, \mu)$ be a measure-preserving transformation of a probability space $(X, \mu)$. Then:
\begin{enumerate}
\item $T$ is ergodic $\iff$ every single (positive- measure) subset family $\mathcal{A} = \{A\}, \mu(A) >0$ is BC $\iff$ every such sequence is sBC.
\item $T$ is weak mixing $\iff$ every finite family $\mathcal{A} = \{A_1, \ldots, A_k\}$ of positive measure sets is BC $\iff$ for every such sequence $S_n/E_n \rightarrow 1$ in the  $L^2$-metric.
\item $T$ is lightly mixing $\iff$ every finite family $\mathcal{A}$ is BC.
\end{enumerate}
We recall that $T$ is \emph{lightly mixing} if for any two sets $A$ and $B$ of positive measure, $\mu(T^{-n}A \cap B)>0$ for all $n>>0$. This condition lies strictly in between mixing and weak mixing.
\end{Theorem}

Notice that both various types of mixing and various types of BC properties can be viewed as different generalizations of ergodicity. However, in general, neither strong mixing nor strong BC properties imply each other, we will see more examples of this in the section on elliptic systems, particularly the results of Fayad~\cite{Fayad}.

However, when we have information on the rate of mixing, we can often derive strong BC properties. For many systems with hyperbolic behavior, we can obtain the estimates on the covariances $|p_{ij} - p_i p_j|$ using exponential estimates on the decay of correlations for the map $T$, and use the Sprindzuk lemma, as in the Kleinbock-Margulis result Theorem~\ref{loglawzero}. This is the method used by Chernov-Kleinbock to prove Theorem~\ref{CKmain}.

Historically, the first result on deterministic Borel-Cantelli properties was Theorem~\ref{Phillips}, due to Phillip~\cite{Phillip}, who used the properties of expanding maps and the Gauss map on the interval, combined with a version of the Sprindzuk lemma in his proof.

The most classical examples of systems with hyperbolic behavior are Anosov diffeomorphisms of compact manifolds and the corresponding symbolic systems, topological Markov chains. Recall that a diffeomorphism $f$ of a Riemannian manifold $M$  is \emph{Anosov}  (or \emph{uniformly hyperbolic}) if there is a $f$-invariant splitting of the tangent bundle $TM = E_u + E_s$ into  \emph{stable} and \emph{unstable} sub-bundles and a constant $\lambda < 1$ such that for all $x \in M$:

\begin{enumerate}

\item $||Df_x(v)|| < \lambda ||v||$ for all $v \in E_s(x)$

\item $||Df^{-1}_x(v)|| < \lambda ||v||$ for all $v \in E_u(x)$

\end{enumerate}

\noindent Here, $E_u(x)$ and $E_s(x)$ denote the fibers of the bundles $E_s$ and $E_u$ over the point $x$. In this context, Dolgopyat proved the following:

\begin{Theorem}\label{Dolgmain}\cite{Dolgopyat} Let $f: M \rightarrow M$ be an Anosov differomorphism preserving a smooth probability measure $\mu$. Then the family of metric balls is strongly Borel-Cantelli, in particular, $f$ has the shrinking target property.\end{Theorem}

Given a transitive matrix $n \times n$ $0-1$ matrix $A$ (i.e., $A^k$ is a positive matrix for some $k >0$), the \emph{topological Markov chain} is given by $$\Sigma = \{ \omega \in \{1, \ldots, n\}^{\Z}: A_{\omega_i \omega_{i+1}} = 1\}$$\noindent (that is, $\Sigma$ consists of doubly infinite sequences with transitions determined by $A$). $\Sigma$ is compact in the product topology, and the left shift $T$ is a homeomorphism. 

The analogue of balls in this setting are \emph{cylinder sets}, which we obtain by fixing the symbols that appear at a finite number of positions. That is, a cylinder $C \subset \Sigma$ is given by an interval $\Lambda = [n, m] \subset \Z$ and a word $\omega^C \in \{1, \ldots, n\}^{m-n+1}$, $$C = \{ \omega \in \Sigma: \omega_{n+i} = \omega^C_i+1 \mbox{ for } 0 \le i \le m-n\}.$$

We can define a metric on $\Sigma$ by fixing an $a \in (0, 1)$ and setting $$d(\omega, \omega^{\prime}) = a^n,$$ where $n = \max{j: \omega_i = \omega^{\prime}_i \forall |i| < j}$. Given any H\"older continuous function (aka \emph{potential}) $\phi: X \rightarrow \R$, we can define an associated shift invariant measure (the \emph{Gibbs measure}) $\mu_{\phi}$. For the explicit construction see~\cite{KH}, for example.

To state the main theorem of Chernov-Kleinbock~\cite{CK}, we require a further definition: we say two intervals $[n_1, m_1]$ and $[n_2, m_2]$ are \emph{$D$-nested} for $D \geq 0$ if either $[n_1, m_1] \subset [n_2 -D, m_2 + D]$ or $[n_2, m_2] \subset [n_1 - D, m_1 +D]$.

\begin{Theorem}\label{CKmain}\cite{CK} Let $\{C_n\}$ be a sequence of cylinder sets defined on intervals $\Lambda_n$. Suppose for all $m, n$, $\Lambda_m$ and $\Lambda_n$ are $D$-nested. Then if $\sum \mu(C_n) = \infty$, $\{C_n\}$ is sBC \end{Theorem}

\section{Elliptic systems}\label{ellip}

At the opposite end of the spectrum from hyperbolic systems are \emph{elliptic} systems, that is, systems where nearby trajectories do not diverge, but rather move in parallel fashion. The most classical examples of elliptic systems are rotations on the circle $S^1$ and torus $\bb{T}^n$. Surprisingly, despite the absence of mixing for these transformations, one can prove many different types of shrinking target results. 

Perhaps less surprisingly, these systems are also very useful for deriving counter-examples for various possible relationships between different types of shrinking target properties and various notions of mixing.

Let $\bb{T}^n = \R^n/ \Z^n$. Given $\alpha \in \R^n$, let $R_{\alpha}$ be the corresponding rotation of $\T^n$. Given $x \in \R^n$, let $||x||_{\Z}$ be the distance from the lattice $\Z^n$. We recall some basic concepts from Diophantine approximation. Given $\sigma \geq 0$, define 

$$\Omega_n(\sigma) = \{ \alpha \in \R^n: \exists C>0 \mbox{ such that } \forall k \in \N ||k\alpha||_{\Z} \geq Ck^{-(1+\sigma)/n}\}$$

\noindent If $\alpha \in \Omega_n(0)$, we say $\alpha$ is a vector of \emph{constant type}. 

The following result was proved by Kurzweil~\cite{Ku}, and rediscovered by Fayad~\cite{Fayad}:

\begin{Theorem}\label{MSTP}~\cite{Ku, Fayad} $R_{\alpha}$ has the MSTP if and only if $\alpha$ is of constant type. \end{Theorem}

Fayad also showed that $R_{\alpha}$ \emph{never} has the STP. In fact, one can produce a sequence of targets with divergent measure such that the $\limsup$ set has measure 0. 

Given $s \geq 1$, we also define the $s$-STP and $s$-MSTP properties as follows; $T$ has the $s$-STP if every sequence of balls $\{B_n\}$ with $\sum \mu(B_n)^s = \infty$ is a Borel-Cantelli sequence  for $\{T^n\}$, i.e., $\mu(x: T^n x \in B_n \mbox{ infinitely often}) = 1$. We say $T$ has the $s$-MSTP if we restrict the requirement to sequences of balls $\{B_n\}$ with decreasing radii. 

Using this terminology, he also shows that $\alpha \notin \Omega_n(1)$ implies that $R_{\alpha}$ does not have the $\frac{n+1}{n}$-MSTP.  

The other main result of~\cite{Fayad} is a construction (using rotations) of a real analytic diffeomorphism of $\bb{T}^3$ that is mixing of all orders, but does not have the MSTP. 

Tseng~\cite{Tseng} studied in particular the case of rotations of the circle $S^1$, proving the following remarkable result, using continued fractions techniques:

\begin{Theorem}\label{sMSTP}~\cite{Tseng} Let $s \geq 1$. Then $R_{\alpha}: S^1 \rightarrow S^1$ has the $s$-MSTP if and only if $\alpha \in \Omega_1(s-1)$.\end{Theorem}

The paper~\cite{Tseng} also contains a proof of the following result, which was also found by Kim~\cite{Kim}:

\begin{Theorem}\label{liminf} For all $\alpha \in \R\backslash\Q$, $$\liminf_{n \rightarrow \infty} n ||n\alpha -s||_{\Z} = 0$$ \noindent for almost every $s \in \R$. \end{Theorem}

As a corollary of this result, one obtains a logarithm law for circle rotations:

\begin{Cor}\label{circleloglaw}  For all $\alpha \in \R\backslash\Q$, for almost all $x \in S^1$, $$\limsup_{n \rightarrow \infty} \frac{-\log||R_{\alpha}^n x||_{\Z}}{\log n} = 1.$$\end{Cor}

\section{Parabolic systems}\label{para}

There are, of course, many systems with intermediate behavior: neither elliptic nor hyperbolic. Some examples of these so-called \emph{parabolic} systems include interval exchange transformations, billiards in Euclidean polygons, and unipotent actions on homogeneous spaces. 

\subsection{Interval Exchange Transformations}

Recall that an \emph{interval exchange transformation} (IET) is given by cutting the interval $[0, 1]$ into $m$ subintervals of lengths given by a vector $\lambda = (\lambda_1, \ldots, \lambda_m)$, and a permutation $\pi \in S_m$. The transformation $T = T_{\lambda, \pi} : [0, 1] \rightarrow [0,1]$ is given by gluing the subintervals together in the order given by $\pi$ and preserving the orientation. These are natural generalizations of circle rotations (which can be viewed as exchanges of two intervals), and are closely related to flows on flat surfaces and billiards in Euclidean polygons.

In what follows, we assume $m \geq 4$, and that our permutation is not a rotation, i.e., $\pi(i+1) \neq \pi(i)  +1 \mbox{(mod }m)$ for some $i$. We say a permutation $\pi \in S_m$ is \emph{irreducible} if for $1 \le n < m$, $$\pi(\{1, \ldots, n\}) \neq \{1, \ldots, n\}.$$ We say an IET is irreducible if its associated permutation is. IETs clearly preserve Lebesgue measure. Unlike circle rotations, for $m \geq 4$, there are minimal (i.e., all non-singular orbits are dense), non-uniquely ergodic IETs~\cite{KeynesNewton}.

It was conjectured by Keane, and proved by Masur~\cite{Masur1} and Veech~\cite{Veech} independently that for an irreducible permutation $\pi$, that $T_{\lambda, \pi}$ is in fact uniquely ergodic   (that is, Lebesgue measure is the \emph{only} preserved measure) for almost all $\lambda \in \Delta^n$, with respect to Lebesgue measure on the simplex $\Delta_n = \{ x \in \R_+^n: \sum x_i = 1\}$. More recently, Avila-Forni~\cite{AF} showed that a generic IET is in fact weak-mixing. 

These proofs of these results use in a crucial fashion various renormalization procedures for IETs, involving induced maps on certain subintervals, and closely related to Teichm\"uller geodesic flow. This discussion is beyond the scope of our paper, but we refer the interested reader, to, e.g.~\cite{Yoccoz}, for an excellent survey. 

The natural question in our context  is: What type of shrinking target property can one prove for IETs? There has been some recent activity in the area; much of it, however is as yet unpublished. Boshernitzan-Chaika recently showed:

\begin{Theorem}\label{chaika}~\cite{BChaika} Let $\pi \in S_m$ be an irreducible permutation, and $\{r_n\}$ a sequence of positive real numbers with $\sum r_n = \infty$. For almost every triple $(\lambda, x, y) \in \Delta^m \times [0, 1]^2$, $$T_{\lambda, \pi}^n x \in B(y, r_n) \mbox{ infinitely often.}$$\end{Theorem}

Note that since the set off full measure \emph{depends} on the sequence, we cannot say that for a generic IET, almost every sequence of targets $B(y, r_n)$ is Borel-Cantelli. 

There is a result in this direction by Galotolo-Kim~\cite{GK}, who restricted their attention to sequences of radii satisfying $r_n = o(1/n)$. Another result in this direction, also due to Boshernitzan-Chaika, is as follows:

\begin{Theorem}\label{BoshChai2}~\cite{BChaika} Let $T = T_{\lambda, \pi}$ be a minimal IET, and $\mu$ an invariant measure. Then for $\mu \times \mu$-a.e. $(x,y) \in [0,1] \times [0, 1]$, $$\liminf_{n \rightarrow \infty} n|T^n x - y| = 0.$$\end{Theorem} 

However, if we ask about a specific sequence of targets (that is, targets centered at a particular point $y$), even this result does not yield any information.

In recent discussions, using the notion of `balanced' times under the Rauzy induction map, the author and C.~Ulcigrai~\cite{AUlcigrai} proved the following analogue of Corollary~\ref{circleloglaw}, giving some information about targets near $0$, and a logarithm law for IETs:

\begin{Theorem}\label{IETloglaw}. Let $\pi \in S_m$ be an irreducible permutation. Then for almost every  $\lambda \in \Delta^m$, for all non-singular $x \in [0, 1]$ $$\limsup_{n \rightarrow \infty} \frac{-\log T_{\lambda, \pi}^ n x }{\log n} = 1.$$\end{Theorem}

\subsection{Unipotent flows}\label{unipotent}

Logarithm laws for unipotent actions on homogeneous spaces were studied by the author and Margulis in~\cite{AM}. We obtained results on general unipotent flows (in the setting of Theorem~\ref{loglawzero}), actions of \emph{horospherical} subgroups, as well as very precise results for unipotent flows on the space of lattices $SL(n,\R)/SL(n, \Z)$, which we describe in the next subsection~\ref{lattices}. We briefly describe the general result on unipotent flows:

We fix notation as in Theorem~\ref{loglawzero}. 

\begin{Theorem}\label{horobounds}~\cite{AM}Let $\{u_t\}_{t \in \R} \subset G$ denote a one-parameter unipotent subgroup. Then there is a $0 < \alpha \le 1$ such that for $\forall y$,  $\mu$-a.e. $x$, $$ \alpha/k \le \limsup_{t \rightarrow \infty} \frac{d(u_t x, y)}{\log t} \le 1/k.$$\end{Theorem}

This result uses a generalization of the Borel-Cantelli lemma to non-independent systems, but the information about mixing for unipotent flows is not sufficient to allow us to apply the Sprindzuk lemma (for example) and derive precise asymptotics. 

For horospherical actions satisfying certain technical conditions, one can derive more precise results. As an example, let $G = SO(n, 1)$, $\Gamma \subset SO(n,1)$ ($n \geq 2$) be a non-uniform lattice, so $N = \h^n/\Gamma$ is a finite-volume, non-compact hyperbolic manifold. Let $SN = SO(n, 1)/\Gamma$ be the unit tangent bundle. Let $H \subset SO(n, 1)$ denote the subgroup such that if $(x, v) \in SN$ ($x \in N$, $v \in S_x N$), then $H(x,v)$ is the unstable manifold through $(x,v)$ for the geodesic flow $g_t$ on $SN$. 

By definition $H$ is horospherical with respect to $\{g_t\}_{t \in \R} \subset SO(n, 1)$, the one-parameter subgroup corresponding to geodesic flow on $SN$. Let $B \subset H^{+}$ be open and $B_t = g_{\log t} B g_{-\log t}$. Let $d(.,.)$ denote the distance on $N$, and $p: SN \rightarrow N$ be the natural projection.

\begin{Cor}\label{hyphoro} For all $x, y \in N$, for almost all $v \in S_x N$ (with respect to Lebesgue measure on the sphere), $$\limsup_{t \rightarrow \infty} \frac{ \sup_{h \in B_t} d(p(h(x,v)), y)}{\log t} = 1.$$ Moreover, for all $(x, v) \in SN$ such that $H(x, v)$ is not periodic,  $$\limsup_{t \rightarrow \infty} \frac{ \sup_{h \in B_t} d(p(h(x,v)), y)}{\log t} \geq 1.$$ \end{Cor}

There is a nice dynamical argument behind this result: consider the piece of horosphere $B_{e^t} x$ we want to show that there are points on it depth $t$ into the cusp. Pulling back under $g_t$, we get the `unit' piece of horosphere $B_1 g_{-t} x$. Since divergent geodesic trajectories are dense, if we 'thicken' this horosphere we obtain nearby points who diverge under $g_t$, i.e., points that will be distance $t$ away at time $t$ (in the hyperbolic setting, geodesics can only diverge at one rate). But these points will be very close to our original horosphere $B_{e^t}x$, since $g_t$ only expands in the horospherical direction, not in the `thickening' direction. This is, of course, an informal argument, and needs to be made precise. Details can be fond in~\cite{AM}.

\subsection{The space of lattices}\label{lattices}

Let $X_n = SL(n, \R)/SL(n, \Z)$ denote the space of unimodular lattices in $\R^n$. Let $\mu = \mu_n$ be Haar measure on $X_n$. Define $\alpha_1 : X_n \rightarrow \R^+$ by $$\alpha_1 (\Lambda) := \sup_{ 0 \neq v \in \Lambda} \frac{1}{||v||}.$$ Let $\{u_t\}_{t \in \R}$ denote a unipotent  one-parameter subgroup of $SL(n, \R)$. 

\begin{Theorem}\label{loglawlattices}~\cite{AM}For $\mu$- a.e. $\Lambda \in X_n$, $$\limsup_{t \rightarrow \infty} \frac{\log \alpha_1(u_t \Lambda)}{\log t} = \frac{1}{n}.$$\end{Theorem}

For $n \geq 3$, this result is based on the following lemma from the geometry of numbers:

\begin{Theorem}\label{measmink}~\cite{AM} Let $n \geq 3$. There is a constant $C_n$ such that  if $A$ is a bounded Borel measurable set in $\R^n$, with $m(A) = a>0$, $$\mu(\Lambda \in X_n: \Lambda \backslash \{0\} \cap A) \le \frac{C_n}{a}$$\end{Theorem}

To derive Theorem~\ref{loglawlattices} from Theorem~\ref{measmink}, one constructs a sequence of measurable subsets $\{A_k\}$ of $\R^n$ with $\mu(A_k) \rightarrow \infty$ so that if a lattice $\Lambda$ has a vector in $A_k$, then there is a time $t_k \approx k$ such that $u_{t_k} \Lambda$ has an appropriately short vector, i.e., a vector of length approximately $k^{-1/n}$.

\subsection{Teichm\"uller horocycle flow}

The author and Yair Minsky proved a similar result for horocycle flow on the space of abelian differentials on surfaces: 

Let \begin{equation}\label{horodef}\left\{h_t =\left(\begin{array}{cc} 1 & t \\ 0 & 1
\end{array}\right)\right\}_{t \in \R}\end{equation}

\begin{Theorem}\label{strata}~\cite{AMinsky}Let $\hh$ denote a stratum of the space of abelian differentials on a surface of genus $g >1$. Let $\lambda: \hh \rightarrow \R^+$ be defined by $$\lambda(\omega) = \sup_{v \in V_{sc}(\omega)} \frac{1}{||v||},$$ where $V_{sc}(\omega)$ is the set of (holonomy vectors of) saddle connections on $\omega$. Then for almost every (with respect to the Lebesgue measure) $\omega \in \hh$, $$\limsup_{t \rightarrow \infty} \frac{\log \lambda(h_t \omega)}{\log t} = 1/2.$$\end{Theorem}

Similar to the proof of Theorem~\ref{loglawlattices}, the proof this result relies on finding saddle connections with holonomy vectors in appropriate regions of the plane $\R^2$.

\section{Conjectures and open problems}\label{conj}

We conclude by outlining a selection of interesting open problems and directions to pursue:

\begin{enumerate}

\item What shrinking target properties can we prove for IETs? What is the precise relationship between them and the integrability properties of the Kontsevich-Zorich cocycle (see, for example, ~\cite{AF} for a precise description of this cocycle, and its relationship to the weak-mixing properties of IETs)?

\item What shrinking target properties can we prove for the actions of products of linear algebraic groups over local fields (i.e., $S$-algebraic groups)?

\item What is the excursion behavior into cusps of stable/unstable manifolds for various (partially) hyperbolic systems? In particular, for the action of horospherical subgroups on $G/\Gamma$, and spaces of measured foliations in Teichmuller spaces.

\item In many of the above results, the shrinking targets are neighborhoods of infinity, i.e., contained in the cusp of $G/\Gamma$ or moduli space. What shrinking target properties can be proved for sets which are contained in the compact part of these spaces (as in Theorems~\ref{mac},~\ref{Dolgmain})?

\item In what non-hyperbolic contexts can Borel-Cantelli properties be extended to Strong Borel-Cantelli properties?

\item A meta-question: in all of the above results, either the base dynamical system or an associated \emph{renormalization} system (for example, for IETs, the renormalization system is Teichm\"uller geodesic flow) exhibits hyperbolic behavior (and in particular, rapid mixing). Does this always need to be the case for a dynamical system to have shrinking target properties?

\end{enumerate}

\end{document}